# Curse-of-dimensionality revisited: Collapse of the particle filter in very large scale systems

Thomas Bengtsson[1], Peter Bickel[2] and Bo Li[*,3]

*Bell Labs, University of California, Berkeley, and Tsinghua University*

**Abstract:** It has been widely realized that Monte Carlo methods (approximation via a sample ensemble) may fail in large scale systems. This work offers some theoretical insight into this phenomenon in the context of the particle filter. We demonstrate that the maximum of the weights associated with the sample ensemble converges to one as both the sample size and the system dimension tends to infinity. Specifically, under fairly weak assumptions, if the ensemble size grows sub-exponentially in the cube root of the system dimension, the convergence holds for a single update step in state-space models with independent and identically distributed kernels. Further, in an important special case, more refined arguments show (and our simulations suggest) that the convergence to unity occurs unless the ensemble grows super-exponentially in the system dimension. The weight singularity is also established in models with more general multivariate likelihoods, e.g. Gaussian and Cauchy. Although presented in the context of atmospheric data assimilation for numerical weather prediction, our results are generally valid for high-dimensional particle filters.

## 1. Introduction

With ever increasing computing power and data storage capabilities, very large scale scientific analyses are feasible and necessary (e.g. [8]). One important application area of high-dimensional data analysis is the atmospheric sciences, where solutions to the general (inverse) problem of combining data and model quantities are commonly required. For instance, to produce real-time weather forecasts (including hurricane and severe weather warnings), satellite radiance observations of humidity and radar backscatter of sea surface winds must be combined with previous numerical forecasts from atmospheric and oceanic models. To such ends, recent work on numerical weather prediction is cast in probabilistic or Bayesian terms [10, 21, 26], and much focus in the literature on the assimilation of data and numerical models pertains to the sampling of high-dimensional probability density functions (pdf) [1, 18, 25, 27]. Motivated by these sampling techniques, we investigate the dangers of naively using Monte Carlo approximations to estimate large

---

[*]Supported in part by National Natural Science Foundation of China (70621061).
[1]Bell Labs, 600 Mountain Avenue, Murray Hill, NJ 07974, USA, e-mail: tocke@research.bell-labs.com
[2]University of California, Department of Statistics, 367 Evans Hall #3860, Berkeley, CA 94720-3860, USA, e-mail: bickel@stat.berkeley.edu
[3]Tsinghua University, 440 Weilun Hall, School of Economics and Management, Beijing 100084, China, e-mail: libo@sem.tsinghua.edu.cn

scale systems. Specifically, in the context of the particle filter, we show that accurate estimation of (truly) high-dimensional pdfs require ensemble sizes that essentially grow exponentially with the system dimension.

Some recent work in numerical weather prediction has extended Kalman filter solutions to work efficiently in Gaussian systems with degrees of freedom exceeding $10^6$. One popular extension is given by the ensemble Kalman filter, a Monte Carlo based filter version which draws samples from the posterior distribution of the atmospheric state given the data and the model [6, 11]. However, the task of sampling in real-time from such high-dimensional systems is conceptually non-trivial: computational resources limit sample sizes to several orders of magnitude smaller than the system dimension. To address sampling errors associated with small ensembles, various approaches leverage sparsity constraints to attenuate spurious correlations [15, 17, 25]. Moreover, in the Gaussian case, for systems with a finite number of dominant modes, moderate sample sizes are sufficient to accurately estimate posterior means and covariances [12].

For longer forecast lead times, the involved dynamical models exhibit strongly non-linear behavior and produce distinctly non-Gaussian forecast distributions (e.g. Figure 2, [3]). In these situations, optimal filtering requires the use of more fully Bayesian filtering methods to combine data and models. In the context of oceanographic data assimilation, one such approach is considered by [27], who proposes a sequential importance sampling algorithm to obtain posterior estimates of oceanic flow structures. This method falls within the set of procedures typically referred to as particle filters (e.g. [9]). Based on a finite set of sample points with associated sample-weights, the particle filter seeks to propagate the probability distribution of the unknown state forward in time using the system dynamics. Once new data is available, Bayes theorem is used to re-normalize the weights based on how "close" the associated sample points are to the data.

Although successfully applied to a variety of settings, particle filters often yield highly varying importance weights and are known to be unstable even in low-order models. Remedies to stabilize the filter include re-sampling (re-normalizing) the involved empirical measure at regular time intervals [14, 19], marginalizing or restricting the sample space [20, 22], and diversifying the sample (e.g. [13]). However, these approaches serve to improve filter performance in low-dimensional systems, but do not fundamentally address slow convergence rates when the particle filter is applied in large scale systems. In particular, as noted by e.g. [27] and [1], when applied to geophysical models of high dimension, sequential importance sampling collapses to a point mass after a few (or even one!) observation cycles. To shed light on the effects of dimensionality on filter stability, our work provides necessary sample size requirements to avoid weight degeneracies in truly large scale problems.

This work is outlined as follows. The next section formulates the problem of using ensemble methods for approximation purposes in large scale systems, and provides motivating simulations illustrating the potential difficulties of high-dimensional estimation. Our main developments are then presented in Section 3 where we give general conditions under which the maximum of the sample weights in the (likelihood based) particle filter converges to one if the ensemble size is sub-exponential in the cube root of the system dimension. The convergence is established in a Gaussian context, but is also argued for observation models with independent and identically distributed (*iid*) kernels. The validity of the weight collapse in the case when the ensemble grows sub-exponentially in the system dimension (as suggested by the simulations) is discussed as an extension to the multivariate Cauchy kernel and its apparent slower collapse. In Section 4, for completeness, in a Gaussian context, we show that the particle filter behaves as desired if the ensemble size is super-exponential in the system dimension. A brief discussion in Section 5 concludes our work.

## 2. Model setting and motivation

This section gives the model setting and describes the Monte Carlo estimation problem under consideration. To illustrate the difficulty of high-dimensional estimation, we also provide motivating examples that describe the weight singularity.

### *2.1. Model setting*

The statistical context in which we motivate our work is as follows. Consider a set of $n$ sample points $\mathbf{X} = \{X_1, \ldots, X_n\}$, where $X_i \in \Re^d$ and both the sample size $n$ and system dimension $d$ are "large". We assume that the sample $\mathbf{X}$ is drawn randomly from the prior (or proposal) distribution $p(X)$. New data $Y$ is related to the state $X$ by the conditional density $p(Y|X)$. For concreteness, a functional relationship $Y = f(X) + \varepsilon$ is assumed, and $\varepsilon$ is taken to be independent the state $X$. The goal is to estimate posterior expectations using the importance ratio: *i.e.*, for some function $h(\cdot)$, we want to estimate

$$E(h(X)|Y) = \int h(X) \frac{p(Y|X)p(X)}{\int p(Y|X)p(X)dX} dX,$$

and use

$$\hat{E}(h(X)|Y) = \sum_{i=1}^{n} h(X_i) \frac{p(Y|X_i)}{\sum_{j=1}^{n} p(Y|X_j)}$$

as an estimator. Based on this formulation, the weights attached to each ensemble member

(1) $$w_i = \frac{p(Y|X_i)}{\sum_{j=1}^{n} p(Y|X_j)}$$

are the primary objects of our study. As mentioned, in large scale analyses, the weights in (1) are highly variable and often produce estimates $\hat{E}(\cdot)$ which are collapsed onto a point mass with $max(w_i) \approx 1$. For high-dimensional systems, this degeneracy is pervasive and appears to hold for a wide variety of prior and likelihood distributions.

Next we illustrate the degeneracy of the sample weights as the dimension of $X$ and $Y$ grows large.

### *2.2. Motivating examples*

To illustrate the weight collapse of the particle filter we simulate weights from (1) for both a Gaussian and a Cauchy distributional setting. These densities were chosen to parallel the work of [27], which attempts to address particle filter collapse by modeling the observation noise using a heavy-tailed distribution.

In our simulations, the $d \times 1$ observation vector $Y$ is related to the unobserved state variable $X$ through the model $Y = HX + \varepsilon$. Here, the observation operator $H$ is set equal to the $d \times d$ identity matrix, denoted $I_d$, and the proposal (i.e. forecast) distribution is taken as zero-mean Gaussian with covariance $cov(X) = I_d$, denoted $X \sim N(0, I_d)$. These choices for $H$ and $cov(X)$ allow us to straightforwardly manipulate the system dimension, and to study the behavior of the maximum weight for various choices of $d$ and $n$. For the Gaussian case we let $\varepsilon \sim N(0, I_d)$, while for the Cauchy case we investigate two scenarios. First, the components of $\varepsilon = [\epsilon_1, \ldots, \epsilon_d]^T$ are taken as *iid* Cauchy, where each component has location and scale parameters set equal to zero and one, respectively. Second, $\varepsilon$ is taken as multivariate Cauchy, again with location and scale parameters equal to zero and one.

In each simulation run of the maximum weight, an observation $Y$ is drawn from $p(Y)$ and a random sample of $n$ particles $X_1, \ldots, X_n$, is obtained from $p(X) = N(0, I_d)$. Then, given $Y$, the sample particles are re-weighted according to (1) and the maximum weight, denoted $w_{(n)}$, is determined. To evaluate the dependence of $w_{(n)}$ on $d$ and $n$, we vary the system dimension at three levels, and let the Monte Carlo sample size increase polynomially in $d$ at a rate of 2.5, *i.e.* we set $n = d^{2.5}$. For the Gaussian setting, we let $d = 10, 50, 100$ and obtain $n = 316, 17676, 100000$, while, for the Cauchy settings, where convergence to unity was observed to be slower, the maximum dimension is increased to $d = 400$. Thus, for the Cauchy simulations, we set $d = 10, 50, 400$, and get $n = 316, 17676, 3200000$. Our setup results in nine simulation sets each for the three distributional settings, and each set is based on 400 independent draws of $w_{(n)}$.

Histograms of the maximum weight $w_{(n)}$ for the Gaussian setting are displayed in Figure 1. The effect of changing the dimension is represented column-wise, and the effect of changing the Monte Carlo sample size is given row-wise. Each plot also depicts the corresponding sample mean maximum weight $\hat{E}(w_{(n)})$ as a vertical line.

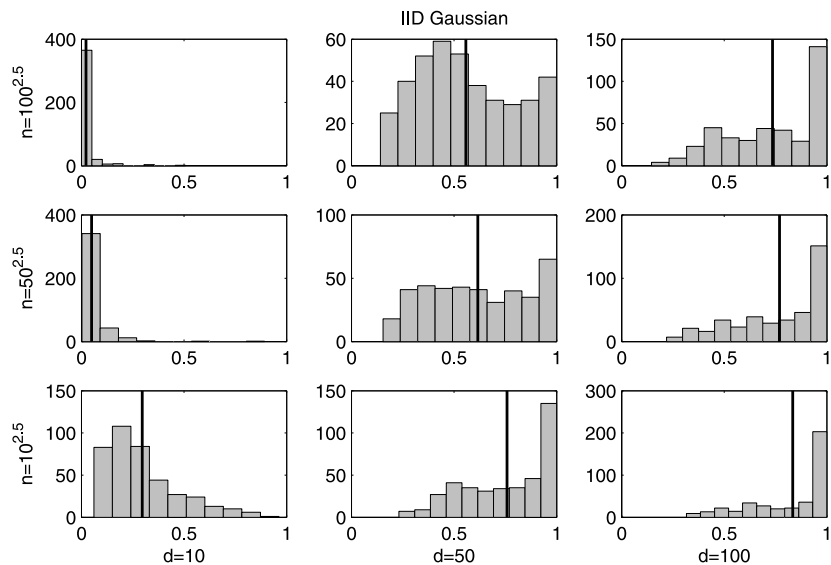

FIG 1. *Sampling distribution of $w_{(n)}$ for the Gaussian case. The nine plots show histograms of $w_{(n)}$ with system dimension varied column-wise ($d = 10, 50, 100$) and sample size varied row-wise ($n = 10^{2.5}, 50^{2.5}, 100^{2.5}$). The vertical line in each plot depicts $\hat{E}(w_{(n)})$. Each histogram is based on 400 simulations.*

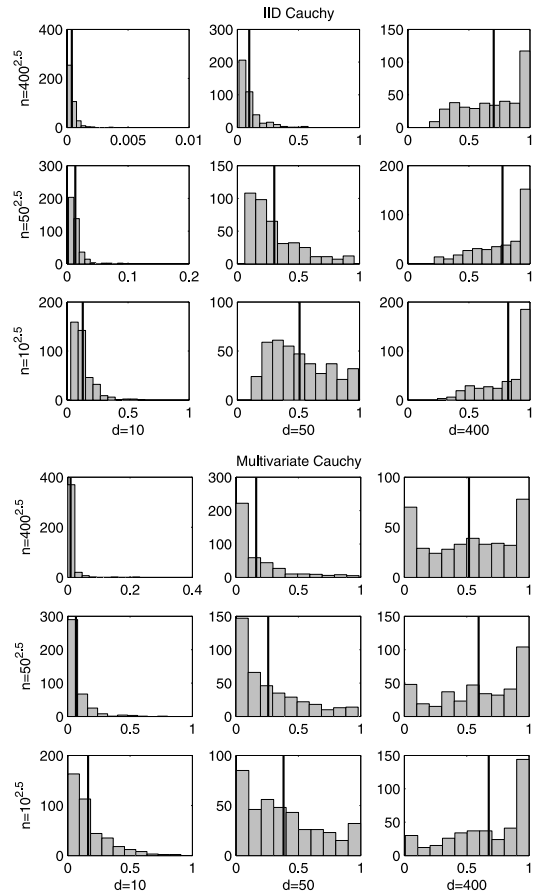

FIG 2. *Sampling distribution of $w_{(n)}$ for $\varepsilon$ iid Cauchy (top panel), and $\varepsilon$ multivariate Cauchy (bottom panel). In each panel, the system dimension d is varied column-wise ($d = 10, 50, 400$), and the sample size is varied row-wise ($n = 10^{2.5}, 50^{2.5}, 400^{2.5}$). The vertical line in each plot depicts $\hat{E}(w_{(n)})$. Each histogram is based on 400 simulations.*

As indicated, for a fixed sample size (i.e. within a row), the distribution of $w_{(n)}$ is dramatically shifted to the right, and we see that $w_{(n)}$ approaches unity as $d$ grows. Moreover, the same singularity is evident along the diagonal from the lower left ($d = 10, n = 316$) to the upper right ($d = 100, n = 100000$) histogram. Hence, even as $n$ grows at a polynomial rate in $d$, $w_{(n)}$ still approaches unity.

The histograms in the two panels of Figure 2 show the simulation results for the Cauchy settings. As in the previous figure, the effect of changing the dimension is given column-wise, and the effect of varying the ensemble size is given row-wise. Also depicted is the sample mean maximum weight $\hat{E}(w_{(n)})$ as a vertical line. For the *iid* Cauchy case, displayed in the top panel, for fixed $n$, the sampling distribution of $w_{(n)}$ is clearly shifted to the right as $d$ increases. Moreover, similarly to the Gaussian results, the weight singularity is again evident in the histograms along the diagonal where the sample size grows as $n = d^{2.5}$. As with the *iid* setting, the histograms of $w_{(n)}$ for the multivariate Cauchy case (bottom panel) also demonstrate weight collapse. However, as can be seen, the convergence is slower. The reasons for the apparent slower collapse will be discussed in Section 3.3.

To shed light on the illustrated weight collapse, we next present a formal study of the behavior of the maximum importance weight. Specifically, the next section develops conditions under which the maximum weight $w_{(n)}$ converges to unity.

## 3. The singularity of the maximum weight for models with *iid* kernels

In this section, we make precise the reasons for the previously described weight collapse. The primary focus is on situations where the likelihood $p(Y|X)$ is based on *iid* components (or *iid* blocks) and the proposal distribution $p(X)$ is Gaussian. Our basic tool, given below in Lemma 3.1, gives reasonable conditions under which the maximum weight based on the general form in (1) converges to unity as $d$ and $n$ grow to infinity.

Our main insight is that, for large $d$, $p(Y|X_i)$ is often well approximated by the form

$$p(Y|X_i) \sim exp\{-(\mu d + \sigma\sqrt{d}Z_i)\}(1 + o_p(1)), \quad (2)$$

where $Z_i$ follows the standard normal distribution and where $\mu$ and $\sigma$ are positive constants. The justification of this form is highlighted by the developments in Sections 3.1 – 3.3, and the conditions under which (2) produces weight collapse are made precise below in Lemma 3.1. Note that, with $Z_{(i)}$ representing the $i$:th order statistic from an ensemble of size $n$, the maximum weight $w_{(n)}$ behaves as

$$w_{(n)} \sim \frac{1}{1 + \sum_{\ell=2}^{n} e^{-\sigma\sqrt{d}(Z_{(\ell)} - Z_{(1)})}}.$$

Thus, to establish weight collapse, it suffices to show that the denominator in the above expression converges to unity for large $d$. The following lemma gives the strong version of (2), which is needed for our conclusion, and formalizes the convergence.

**Lemma 3.1.** *Let $S_i, i = 1, \ldots, n$, be independent random variables with cumulative distribution function (cdf) $G_d(\cdot)$ satisfying*

$$G_d(s) = (1 + o(1))\Phi(s)$$

*over a suitable interval, where $\Phi(\cdot)$ is the cdf of a standard normal distribution. Specifically, we assume there exist two sequences $\alpha_{n,d}^l$ and $\alpha_{n,d}^u$, where both $\alpha_{n,d}^l$ and $\alpha_{n,d}^u$ tend to 0 as $n$ and $d$ go to infinity, such that, for $s \in [-d^\eta, d^\eta]$ with $0 < \eta < 1/6$, we have*

$$(1 + \alpha_{n,d}^l)\Phi(s) \leq G_d(s) \leq (1 + \alpha_{n,d}^u)\Phi(s). \quad (3)$$

*Let $S_{(1)} \leq \cdots \leq S_{(n)}$ be the ordered sequence of $S_n, \ldots, S_n$, and define, for some $\sigma > 0$, $T_{n,d} = \sum_{\ell=2}^{n} e^{-\sigma\sqrt{d}(S_{(\ell)} - S_{(1)})}$. Then, if $\frac{\log n}{d^{2\eta}} \to 0$,*

$$E(T_{n,d}|S_{(1)}) = o_p(1).$$

A proof of the lemma is provided in the Appendix.

An immediate implication of this result is weight collapse, i.e. if $\frac{\log n}{d^{2\eta}} \to 0$, we have $w_{(n)} \xrightarrow{P} 1$.

With the normality condition (3) valid for $\eta < \frac{1}{2}$, our conclusion would hold, effectively, whenever $\frac{\log n}{d} \to 0$; rather than, as shown above, for $\frac{\log n}{d^{1/3}} \to 0$. Unfortunately, unless $G_d(s) \equiv \Phi(s)$, our proof is valid only for $\eta < \frac{1}{6}$. However, using more refined arguments in conjunction with Lemma 2.5 of [23] one can show that, under the conditions of Lemma A.1 (see Appendix), only $\frac{\log n}{d^{1/2}} \to 0$ is required for collapse. Furthermore, using specialized arguments, we can show that in the Gaussian prior-Gaussian likelihood setting, $\frac{\log n}{d} \to 0$ is the necessary condition for collapse. We do not describe these arguments in detail, but focus instead on Lemma 3.1 which can be applied more broadly.

We next turn our attention to high-dimensional, linear Gaussian systems.

### 3.1. Gaussian case

We assume a data model given by $Y = HX + \varepsilon$, where $Y$ is a $d \times 1$ vector, $H$ is a known $d \times q$ matrix, and $X$ is a $q \times 1$ vector. Both the proposal distribution and the error distribution are Gaussian with $p(X) = N(\mu_X, \Sigma_X)$ and $p(\varepsilon) = N(0, \Sigma_\varepsilon)$, and the noise $\varepsilon$ is taken independent of the state $X$. Since the data model can be pre-rotated by $\Sigma_\varepsilon^{-1/2}$, we set $\Sigma_\varepsilon = I_d$ without loss of generality (wlog). Moreover, since $EY = EHX$, we can replace $X_i$ by $(X_i - EX_i)$ and $Y$ by $(Y - EY)$ and leave $p(Y|X)$ unchanged. Hence, wlog we also set $\mu_X = 0$. Moreover, define, for conformable $A$ and $B$, the inner product $\langle A, B \rangle = A^T B$ (where the superscript $T$ denotes matrix transpose), and let $\|A\|^2 = \langle A, A \rangle$.

With $p(Y|X) \sim N(HX, I_d)$, the weights in (1) can be expressed as

$$w_i = \frac{exp\big(-\|Y - HX_i\|^2/2\big)}{\sum_{j=1}^n exp\big(-\|Y - HX_j\|^2/2\big)}. \tag{4}$$

To establish weight collapse for high-dimensional Gaussian $p(Y|X)$ and $p(X)$, we show that, under reasonable assumptions, the exponent in (4) satisfies the conditions of Lemma 3.1.

Let $d' = rank(H)$. With $\lambda_1^2, \ldots, \lambda_{d'}^2$ the singular values of $cov(HX)$, define the $d' \times d'$ matrix $D = diag(\lambda_1, \ldots, \lambda_{d'})$. Then, with $Q$ an orthogonal matrix obtained by the spectral decomposition of $cov(HX)$, define the $d' \times 1$ vector $V$ such that

$$Q^T HX = DV.$$

Note that $V_i$ corresponding to $X_i$ is $N(0, I_{d'})$. Since $Q$ is orthogonal, we can write

$$\|Y - HX_i\|^2 = \|Q^T Y - DV_i\|^2 = \sum_{j=1}^{d'} \lambda_j^2 W_{ij}^2 + \sum_{j=d'+1}^{d} \epsilon_{0j}^2, \tag{5}$$

where, conditional on $Y$, $[W_{i1}, \ldots, W_{id'}]^T$ is $N(\xi, I_{d'})$, and where $\epsilon_{0j}$ is the $j$:th component of the observation noise vector $\varepsilon$. The mean vector $\xi = [\mu_1, \ldots, \mu_{d'}]^T$ is given by

$$\xi = D^{-1} Q^T Y = V + D^{-1} \varepsilon', \tag{6}$$

where $V$ and $\varepsilon'$ are independent $N(0, I_{d'})$.

Now, for $i = 1, \ldots, n$, define

$$S_i = \frac{\sum_{j=1}^{d'} \lambda_j^2 (W_{ij}^2 - (1 + \mu_j^2))}{\big(2 \sum_{j=1}^{d'} \lambda_j^4 (1 + 2\mu_j^2)\big)^{1/2}}. \tag{7}$$

Note that the second term in (5) is constant for every $X_i$, and will not influence the weight $w_i$.

We assume,

A1. There is a positive constant $\delta$ such that $\delta < \lambda_1, \cdots, \lambda_{d'} \leq \frac{1}{\delta}$; and

A2. $\sigma_{d'}^2 = \frac{2}{d'} \sum_{j=1}^{d'} \lambda_j^4 (1 + 2\mu_j^2) \to \sigma^2 > 0$.

With these assumptions, Lemma A.2 of the Appendix establishes that the uniform normality condition given by (3) holds for the standardized terms defined in (7). The result is valid for any $0 < \eta < 1/6$ and is based on Theorem 2.5 of [23].

Note that, from (6) and (7), we can write

$$\|Y - HX_i\|^2 \propto \sigma\sqrt{d'} S_i(1 + o(1)) + \sum_{j=1}^{d'} \lambda_j^2 (1 + \mu_j^2),$$

where $o(1)$ is independent of the subscript $i$.

We can now state the following proposition.

**Proposition 3.2.** *Under assumptions* A1 *and* A2, *if* $\frac{\log n}{d'^{2\eta}} \to 0$ *with* $0 < \eta < 1/6$, *we have* $w_{(n)} \xrightarrow{P} 1$.

Proposition 3.2 follows by Lemma A.2 (Appendix) and Lemma 3.1.

The above result implies that, unless $n$ grows super-exponentially in $d'^{1/3}$, we have weight collapse. However, as we show in Section 4, the weight singularity is avoided when $d'$ is bounded, or, more generally, when $\log n/d' \to \infty$. To establish the exact boundary of collapse in the Gaussian setting, a closer analysis of the chi-squared distribution (*c.f.* $\|Y - HX_i\|^2$) with $d'$ degrees of freedom is needed. In particular, using the Poisson sum formula for the tails of the gamma distribution, it can be argued that collapse occurs if $\log n/d' \to 0$. Essentially, with $G_d(s) = \Phi(s)$ and for suitable $\sigma^2$, we can show

$$E(T_{n,d}|S_{(1)}) \sim \sqrt{\frac{2\log n}{\sigma^2 d'}},$$

in probability. As mentioned, we do not describe the specifics of this argument, but focus instead on the more general result of Lemma 3.1.

Next we turn our attention to settings where both the observation and state vectors consist of *iid* components.

### *3.2. General* **iid** *kernels*

It may be speculated that the weight singularity can be ameliorated by the use of a heavy-tailed kernel. However, we argue that, as long as the components of the observation noise are *iid* and $H$ is the identity matrix, we still expect weight collapse for large $d$.

The model setting is $Y = X + \varepsilon$. Let $X_{ik}$ be the $k$:th component of $X_i$. We take $X_{ik}$ *iid* with common density $g(\cdot)$, and take the components of $\varepsilon = [\epsilon_1, \ldots, \epsilon_d]^T$ *iid* with common density $f(\cdot)$. Then, with $\psi(\cdot) = \log f(\cdot)$, given $Y = [y_1, \ldots, y_d]^T$, we can write the weights from (1) as

$$w_i = \frac{exp\bigl(\sum_{j=1}^d \psi(y_j - X_{ij})\bigr)}{\sum_{k=1}^n exp\bigl(\sum_{j=1}^d \psi(y_j - X_{kj})\bigr)}.$$

Define $V_{ij} = \psi(y_j - X_{ij})$, and let $\mu(y_j) = E(V_{ij}|y_j)$ and $\sigma^2(y_j) = E(V_{ij}^2|y_j) - \mu^2(y_j)$, where the expectations are evaluated under $f(\cdot)$. With these quantities, let

$$S_i = \frac{\sum_{j=1}^d (V_{ij} - \mu(y_j))}{\bigl(\sum_{j=1}^d \sigma^2(y_j)\bigr)^{1/2}},$$

for $i = 1, \ldots, n$. Again, let $S_{(1)} \leq \cdots \leq S_{(n)}$ be the ordered sequence of $S_1, S_2, \ldots, S_n$. Now, in analogy to Proposition 3.2: if $\frac{\log n}{d^{2\eta}} \to 0$, and

(i) the sequence $\{V_{ij} - \mu(y_j)\}_{j=1}^d$ satisfies the conditions of Theorem 2.5 of [23]; and
(ii) $\frac{1}{d}\sum_{j=1}^d \sigma^2(y_j) \xrightarrow{P} \sigma^2 > 0$,

then the maximum weight $w_{(n)}$ converges in probability to 1.

To verify the normality approximation in (2) for $S_i$, it is easy to check that $\sigma^2 = E\sigma^2(y_1) < \infty$, and $\max_j |\psi(y_j - X_{ij})| = o_p(d^{1/2})$ uniformly in $i$. The next proposition gives checkable, albeit strong conditions, leading to weight collapse.

**Proposition 3.3.** *Let the components of $X_0 = [X_{01}, \ldots, X_{0d}]^T$ be iid with density $g(\cdot)$, and let the components of $\varepsilon = [\varepsilon_1, \ldots, \varepsilon_d]^T$ be iid with density $f(\cdot)$. Set $Y = X_0 + \varepsilon$. Let $X_1, \ldots, X_n$ be iid vectors, each with iid components $X_{ik}$ with common density $g(\cdot)$. Assume that $f, g$ are such that*

$$E[f^t(Y_j - X_{1j})] < \infty,$$

*for $|t| < \delta$ with $\delta > 0$. Then, with $T_{n,d} = \sum_{\ell=2}^n e^{-\sigma\sqrt{d}(S_{(\ell)} - S_{(1)})}$, if $\frac{\log n}{d^{2\eta}} \to 0$ for $0 < \eta < 1/6$, we have,*

$$E(T_{n,d}|S_{(1)}) \xrightarrow{P} 0.$$

The result follows from Lemma 3.1 and Lemma A.1. Note that most common $f, g$, e.g. the Gaussian and Cauchy, satisfy the conditions of Proposition 3.3.

If $H$ is not the identity, our conclusion may still hold for $\frac{\log n}{d^{1/3}} \to 0$. For general $H$, set $U_i = HX_i$ and let $U_{ik}$ be the $k$-th component of $U_i$. With this quantity, provided the regularity conditions are satisfied, we may apply Theorem 3.23 of [23] to the terms

$$S_{\psi,i} = \frac{\sum_{j=1}^d \psi(y_j - U_{ij}) - \mu(y_j)}{\sigma\sqrt{d}},$$

where $\mu(y_j)$ and $\sigma$ are suitable constants.

Next we turn our attention to the case when the observation noise vector follows the multivariate Cauchy distribution.

### 3.3. Multivariate Cauchy case

Our developments so far have considered settings where $p(Y|X)$ is of multiplicative form; in particular, model settings with *iid* likelihood kernels (or *iid* blocks of observations). We now discuss extensions of our results to include multivariate likelihood functions. Except for the multivariate Gaussian case, however, which can be addressed by rotation (see Section 3.1), we note that no general result exists that addresses a wide range of multivariate likelihood functions. Here, we focus on the multivariate Cauchy distribution, which, in the context of oceanographic data assimilation, was proposed by [27] to avoid weight collapse.

We still entertain the data model $Y = HX + \varepsilon$ and, as in the previous section, restrict $H = I_d$. Here, we let $p(X) = N(0, I_d)$ and take the noise vector $\varepsilon$ to follow the multivariate Cauchy distribution. We note that the multivariate Cauchy distribution is equivalent to the multivariate t-distribution with 1df (e.g. [2], page 55). Then, given data $Y$ and a sample $X_i \sim N(0, I_d)$, $i = 1, \ldots, n$, the multivariate Cauchy weights are given by

$$w_i = \frac{(1 + \|Y - X_i\|^2)^{-\frac{d+1}{2}}}{\sum_{j=1}^n (1 + \|Y - X_j\|^2)^{-\frac{d+1}{2}}}.$$

Now, with $\|Y-X_{[1]}\| \leq \ldots \leq \|Y-X_{[n]}\|$ the order statistics, we express the weights in familiar form, i.e.

$$w_{(n)} = \left[1 + \sum_{j=2}^{n} exp\left(-\frac{d+1}{2}\left(\log(1+\|Y-X_{[j]}\|^2) - \log(1+\|Y-X_{[1]}\|^2)\right)\right)\right]^{-1}.$$

To argue $w_{(n)} \approx 1$ for large $d$, we note first that the scenario considered here is closely related to the Gaussian prior-Gaussian likelihood case.

**Proposition 3.4.** *Let $X$ be zero-mean Gaussian with $cov(X) = I_d$ and take $\varepsilon$ to follow the multivariate Cauchy distribution. Then, with $Y = X + \varepsilon$, as $d \to \infty$, we get*

$$p(X|Y) \xrightarrow{L} N\left(Y\frac{d}{\|Y\|^2}, (1-\frac{d}{\|Y\|^2})I_d\right).$$

The convergence in Proposition 3.4 is in the sense that the finite dimensional distributions on both sides converge to the same limit. Thus, we reach the somewhat surprising conclusion that the posterior distribution of $X$ given $Y$ is Gaussian, but with parameters depending on the data. The result is proved in the Appendix.

In the Appendix we also detail the heuristic argument that, as in the Gaussian-Gaussian case, the non-central $\chi^2_d(\sigma^2 d)$ distribution behaves sufficiently like $N((\sigma^2+1)d, 2(1+2\sigma^2)d)$ to permit exact replacement in our developments. We then reach the conclusion of slower weight collapse for large $n$ and $d$, as substantiated by our simulations. Specifically, we argue that the average rate needed for collapse is $\sqrt{\frac{\log n}{d}}\left|\log\left(\sqrt{\frac{\log n}{d}}\right)\right| \to 0$.

For the Gaussian setting, the next section shows that weight collapse can be avoided for sufficiently large ensembles

## 4. Consistency of Gaussian particle filter

As a complement to the developments in the previous section, we provide a consistency argument for the type of estimators of $E(h(X)|Y)$ that are under consideration here. The developments are made in a Gaussian context and we consider settings where both $d$ and $n$ are large. Specifically, if $\log n/d \to \infty$, we show consistency of $\sum_{i=1}^n w_i h(X_i)$ as an estimator of $E(h(X)|Y)$.

Suppose we have a random sample $\{X_0, X_1, \ldots, X_n\}$, where $X_i \sim N(0, I_d)$. As in the simulation section, let the data $Y_0$ be collected through the model $Y_0 = X_0 + \varepsilon$, where $\varepsilon \sim N(0, I_d)$ is independent of $X_i$. With $p(Y|X) = N(0, I_d)$, let $\{w_1, \ldots, w_n\}$ be the weights obtained by (1).

Now, choose $X_i^*$ from $\{X_1, \ldots, X_n\}$ with probabilities proportional to $\{w_1, \ldots, w_n\}$. Then, with $\delta(\cdot)$ representing the delta function and $\tilde{p}(X|Y_0) = \sum_{j=1}^n w_j \delta(X_j)$, we have $X_i^* \sim \tilde{p}(X|Y_0)$. Further, the expectation of $h(X^*)$ (under the empirical measure) is given by $E^*(h(X^*)|Y_0) = \sum_{j=1}^n w_j h(X_j)$. With this setup, the following result establishes consistency of $\sum_{i=1}^n w_i h(X_i)$ as an estimator of $E(h(X)|Y)$.

**Proposition 4.1.** *Let $h(\cdot)$ be a bounded function from $\Re^d$ to $\Re$. Define $E^*h(X^*) \equiv \sum_{j=1}^n w_j h(X_j)$, the expectation under the (previously defined) empirical measure $\tilde{p}(X|Y_0)$, and let $E_1(\cdot)$ denote expectation evaluated under the (true) posterior $p(X|Y_0)$. Then, if $\log n/d \to \infty$,*

$$|E^*h(X^*) - E_1 h(X)| \xrightarrow{P} 0.$$

The result is proved in the Appendix. We note that Proposition 4.1 is valid for finite dimensional measures. The next corollary follows as an immediate consequence of the established convergence.

**Corollary 4.2.** *Let $\{X_1^*, \ldots, X_n^*\}$ be a random sample from $p^*(X|Y_0)$. Further, with $k$ fixed, let $\nu(x_1^*, \ldots, x_k^*)$ be any random variable depending on the first $k$ coordinates of $X_i^*$. Then, with $\delta_\nu(\cdot)$ representing the delta function for $\nu(\cdot)$, as $\log n/d \to \infty$,*

$$\frac{1}{n}\sum_{j=1}^n \delta_\nu(x_1^*, \ldots, x_k^*) \xrightarrow{L} N\big(\frac{1}{2}[y_{01}, \ldots, y_{0k}]^T, \frac{1}{2}I_k\big).$$

The results extend to include $Y = HX + \varepsilon$, where $rank(H) = d' < \infty$. Of course, $p(Y|X)$ and $p(Y)$ have to be changed in the developments to accommodate arbitrary $H$, but we again have no collapse provided $\log n/d' \to \infty$.

## 5. Discussion

The collapse of the weights to a point mass (with $max(w_i) \approx 1$) leads to disastrous behavior of the particle filter. One intuition about such weight-collapses is well known, but here made precise in terms of $d$ and $n$: Monte Carlo does not work if we wish to compute $d$-dimensional integrals with respect to product measures. The reason is that we are in a situation where the proposal distribution $p(X)$ and the desired sampling distribution are approximately mutually singular and (essentially) have disjoint support. As a consequence, the density of the desired distribution at all points of the proposed ensemble is small, but a vanishing fraction of density values predominate in relation to the others.

Our developments demonstrate that brute-force-only implementations of the particle filter to describe high-dimensional posterior probability distributions will fail. Our work makes explicit the rates at which sample sizes must grow (with respect to system dimension) to avoid singularities and degeneracies. In particular, we give necessary bounds on $n$ to avoid convergence to unity of the maximum importance weight; and, naturally, accurate estimation will require even larger sample sizes than those implied by our results. Not surprisingly, weight degeneracies have been observed in geophysical systems of moderate dimension ([1, 3]; also, C.Snyder/NCAR & T. Hamill/NOAA, personal communication, 2001). The usual manifestation of this degeneracy are Monte Carlo samples that are too "close" to the data, quickly producing singular probability measures, in particular as the filter is cycled over time.

The obvious remedy to this phenomenon is to achieve some form of dimensionality reduction, and the high-dimensional form in which the data are presented is typically open to such reduction with subsequent effective analysis. For instance, in the case of the ensemble Kalman filter, by imposing sparsity constraints through spatial localization ([15, 17]; see also, [12]). Be that as it may, as shown in this work, for fully Bayesian filter analyses of high-dimensional systems, such reduction becomes essential lest spurious sample variability is to dominate the posterior distribution.

In the context of numerical weather prediction, one approach to dimension reduction may be to condition sample draws on a larger information set. One idea is given by [4], who constructs proposal distributions by incorporating dynamic information in a low-order model. Other examples of geophysically constrained sampling schemes are given by Bayesian Hierarchical Models (e.g. [16, 28]), but require computationally heavy, chain-based sampling and thus do not extend in any obvious manner to real-time applications. A more viable possibility for real-time applications in a large scale system is to break the system into lower-dimensional sets, and then sequentially perform the sampling as in [3]. Another approach may be to condition the draws on lower-dimensional sufficient statistics ([5, 24]). A novel idea to improving convergence in the Baysian filter context is considered by [29] in the context of image matching and retrieval. For the purpose of validating the exclusion of parts of the sample space which appear uninteresting given the data, and to speed up the algorithm, they use information theory to restrict the sample space by explicitly incorporating (drawing) samples of *low* probability.

**Appendix**

We first introduce two lemmas that pertain to uniform normal approximations of the distribution of independent sums. Such sums appear in the formulation of the filter weights throughout our work. Valid for moderately large deviations, the first result (Lemma A.1) is a combination of Theorem 2.5 and Theorem 1.31 in [23] and is stated here without proof. The next result (Lemma A.2) is given for the purpose of verifying the Lyapunov quotients conditions appearing in Lemma A.1.

**Lemma A.1.** *Let $\xi_1, \ldots, \xi_d$ be independent random variables with $E\xi_j = 0$ and $\sigma_j^2 = Var(\xi_j^2)$. Set*

$$S_d = \frac{1}{B_d}(\xi_1 + \cdots + \xi_d),$$

*where $B_d^2 = \sum_{j=1}^d \sigma_j^2$, and define the Lyapunov quotients*

$$L_{k,d} = \frac{1}{B_d^k} \sum_{j=1}^d E|\xi_j|^k, \quad k = 1, 2, \ldots.$$

*If $B_d = \sigma d^{1/2}(1 + o(1))$, $\tau_d = cB_d$ (some $\sigma, c > 0$), and $L_{k,d} \leq k!/\tau_d^{k-2}$, $k = 3, 4, \ldots$, then the cdf of $S_d$, denoted $G_d(\cdot)$, satisfies (some $C > 0$)*

$$\left|\frac{G_d(x)}{\Phi(x)} - 1\right| \leq C|x|^3 d^{1/2}, \quad -d^\eta \leq x \leq 0, \quad \eta < 1/6;$$

*and the survival function, $\bar{G}_d(\cdot) = 1 - G_d(\cdot)$, satisfies*

$$\left|\frac{\bar{G}_d(x)}{\bar{\Phi}(x)} - 1\right| \leq C|x|^3 d^{1/2}, \quad 0 \leq x \leq d^\eta, \quad \eta < 1/6.$$

*Thus, under the outlined conditions, the uniform normal approximations*

$$G_d(s) = (1 + o(1))\Phi(s), \text{ and } \bar{G}_d(s) = (1 + o(1))\bar{\Phi}(s)$$

*hold, respectively, over the intervals $[-\tau_d^\eta, 0]$ and $[0, \tau_d^\eta]$ for $\eta < 1/3$.*

The lemma is the basis for the normality conditions of Lemma 3.1. We note that the condition on $L_{k,d}$ is satisfied in the *iid* case if $Ee^{t\xi_j} < \infty$, for $|t| \leq \delta, \delta > 0$.

**Lemma A.2.** *Let $Z_j, j = 1, \ldots, d$, be iid N(0,1). Let $0 < \delta \leq \lambda_1 \leq \cdots \leq \lambda_d \leq \frac{1}{\delta}$ and $\mu_1, \ldots, \mu_d$ be such that*

(8) $$0 < m < \sigma_d^2 = \frac{2}{d} \sum_{j=1}^d \lambda_j^4 (1 + 2\mu_j^2) \leq M < \infty,$$

*for some constants $m, M$. Then, with $\sigma_d^2$ as defined in (8),*

$$\text{(9)} \quad \sigma_d^{-k/2} \sum_{j=1}^{d} \lambda_j^{2k} E \big|(Z_j + \mu_j)^2 - (1 + \mu_j^2)\big|^k \leq C^k k! d^{-\frac{(k-2)}{2}},$$

*for all $k = 3, 4, \ldots$ and a universal constant $C = C(\delta, m, M)$.*

The remainder of the Appendix is devoted to the proofs of the presented lemmas and propositions. The proofs are given in the order in which the results appear in the text.

### *Proof of Lemma 3.1*

Let $S_j$ $(j = 1, \ldots, n)$ be as defined in the lemma. Before proceeding, we note some important results/facts. First, it is well known that

$$\text{(10)} \quad \sqrt{2 \log n} \Big( S_{(1)} + \sqrt{2 \log n} - \frac{\log \log n + \log(4\pi)}{2\sqrt{2 \log n}} \Big) \xrightarrow{L} U,$$

where $U$ has a Gumbel distribution (see, e.g., [7], p. 475). Second, by assumption, $\frac{\sqrt{2 \log n}}{d^\eta} \to 0$ for $0 < \eta < 1/6$, which allows us to use the normal approximation defined in (3). Third, we make frequent use of Mill's ratio: i.e., $\bar{\Phi}(x) \sim \frac{\phi(x)}{x}$, as $x \to +\infty$.

Now, note that

$$\text{(11)} \quad E(T_{n,d}|S_{(1)}) = \frac{(n-1) \int_{S_{(1)}}^{\infty} exp\big(-\sigma\sqrt{d}(z - S_{(1)})\big) dG_d(z)}{\bar{G}_d(S_{(1)})},$$

since, given $S_{(1)}$, the remaining $(n-1)$ observations are *iid* with cumulative density function (cdf) equal to $G_d(z)/\bar{G}_d(S_{(1)})$, $z \geq S_{(1)}$. By (10) and (3) we see that the denominator of the right hand side converges to 1 in probability. To evaluate the numerator of (11), we break the integral into two parts: the first part yields the integral from $S_{(1)}$ to $d^\eta$, and the second part yields the tail integral from $d^\eta$ to $\infty$. We now show that both integrals converge to 0 in probability.

For the first part, using integration by parts (twice) along with the approximation in (3), we can write

$$\int_{S_{(1)}}^{d^\eta} exp\big(-\sigma\sqrt{d}(z - S_{(1)})\big) dG_d(z)$$
$$= G_d(d^\eta) exp\big(-\sigma\sqrt{d}(d^\eta - S_{(1)})\big) - G_d(S_{(1)})$$
$$\quad + \int_{S_{(1)}}^{d^\eta} exp\big(-\sigma\sqrt{d}(z - S_{(1)})\big) G_d(z) dz$$
$$= \big(G_d(d^\eta) - \Phi(d^\eta)\big) exp\big(-\sigma\sqrt{d}(d^\eta - S_{(1)})\big) - \big(G_d(S_{(1)}) - \Phi(S_{(1)})\big)$$
$$\quad + (1 + o(1)) \int_{S_{(1)}}^{d^\eta} exp\big(-\sigma\sqrt{d}(z - S_{(1)})\big) \phi(z) dz.$$

Now, again by assumption,

$$\text{(12)} \quad (n-1) G_d(d^\eta) exp\big(-\sigma\sqrt{d}(d^\eta - S_{(1)})\big) \leq (n-1) G_d(d^\eta) \sim (n-1) \Phi(d^\eta) = o(1).$$

Analogously, $(n-1)\Phi(d^\eta)exp\big(-\sigma\sqrt{d}(d^\eta - S_{(1)})\big) = o_p(1)$. Hence, we have derived that

(13) $\qquad (n-1)\big(G_d(d^\eta) - \Phi(d^\eta)\big)exp\big(-\sigma\sqrt{d}(d^\eta - S_{(1)})\big) = o_p(1).$

Further, by (10) and (3),

(14) $\qquad (n-1)\big(G_d(S_{(1)}) - \Phi(S_{(1)})\big) = o(1)(n-1)G_d(S_{(1)}).$

Let $S$ be an *iid* copy of $S_1$. Then $G_d(S)$ is uniformly distributed on [0,1], and

(15) $\quad E\big(G_d(S_{(1)})\big) = E\big(P(S \leq S_1, \ldots, S \leq S_n|S)\big) = E\big(([\bar{G}_d(S)]^n\big) = \dfrac{1}{n+1}.$

Combining (14) and (15) yields

(16) $\qquad (n-1)\big(G_d(S_{(1)}) - \Phi(S_{(1)})\big) = o_p(1).$

Moreover, by assumption $\sigma\sqrt{d} + S_{(1)} \sim \sigma\sqrt{d} - \sqrt{2\log n} \to \infty$. Thus, we get, again in light of (10),

$$\begin{aligned}
& (n-1)\int_{S_{(1)}}^{d^\eta} exp\big(-\sigma\sqrt{d}(z - S_{(1)})\big)\phi(z)dz \\
\leq\ & (n-1)\,exp(\sigma\sqrt{d}S_{(1)} + \sigma^2 d/2)\bar{\Phi}(S_{(1)} + \sigma\sqrt{d}) + \sqrt{d}(n-1)\bar{\Phi}(d^\eta) \\
\sim\ & \dfrac{(n-1)exp(\sigma\sqrt{d}S_{(1)} + \sigma^2 d/2)\phi(S_{(1)} + \sigma\sqrt{d})}{S_{(1)} + \sigma\sqrt{d}} + \sqrt{d}(n-1)\bar{\Phi}(d^\eta) \\
=\ & \dfrac{(n-1)\,exp\big(-S_{(1)}^2/2\big)}{\sqrt{2\pi}(S_{(1)} + \sigma\sqrt{d})} + \sqrt{d}(n-1)\bar{\Phi}(d^\eta) \\
=\ & o_p(1),
\end{aligned}$$

The implication of the last inequality, in conjunction with (13) and (16), yields

$$(n-1)\int_{S_{(1)}}^{d^\eta} exp\big(-\sigma\sqrt{d}(z - S_{(1)})\big)dG_d(z) = o_p(1).$$

By assumption, the second part of the integral in (11) is bounded by

$$\sigma\sqrt{d}(n-1)\bar{G}_d(d^\eta) = \sigma\sqrt{d}(n-1)\bar{\Phi}(d^\eta) \sim \dfrac{\sigma\sqrt{d}(n-1)\phi(d^\eta)}{d^\eta} = o(1).$$

We have shown that the numerator of the right hand side of (11) converges to 0 in probability. This completes the proof.

**Remark.** If the normal approximation is good enough to avoid the left boundary term in the integration by parts, we believe the convergence rate to be dominated by the quantity

$$\int_{S_{(1)}}^\infty exp\big(-\sigma\sqrt{d}(z - S_{(1)})\big)\,d\Phi(z).$$

As can be seen in the proof, by (10),

$$\int_{S_{(1)}}^\infty exp\big(-\sigma\sqrt{d}(z - S_{(1)})\big)\,d\Phi(z) \sim \dfrac{C(n-1)exp\big(-S_{(1)}^2/2\big)}{\sqrt{2\pi}(S_{(1)} + \sigma\sqrt{d})} = O_p\bigg(\sqrt{\dfrac{\log n}{d}}\bigg).$$

*Proof of Proposition 3.4.* Let $Z_i \overset{iid}{\sim} N(0,1), i = 0, \ldots, d$, and define $\varepsilon = [Z_1, \ldots, Z_d]^T/|Z_0|$. That is, $\varepsilon$ follows the multivariate Cauchy distribution. Let the data be defined by $Y = X + \varepsilon$.

Conditional on $(Y, Z_0)$, the posterior distribution of $X$ is Gaussian, i.e.

$$p(X|Y, Z_0) = N\left(\frac{|Z_0|^2}{1+|Z_0|^2}Y, \frac{1}{1+|Z_0|^2}I_d\right).$$

We approximate

$$\begin{aligned}
\|Y\|^2 &= \|X\|^2 + 2\langle X, \varepsilon\rangle + \|\varepsilon\|^2 \\
&= d + O(\sqrt{d}) + 2O(\sqrt{d}) + \frac{d + O(\sqrt{d})}{|Z_0|^2}, \\
&= d(1 + O(1/\sqrt{d}) + \frac{1 + O(1/\sqrt{d})}{|Z_0|^2}),
\end{aligned}$$

We have $\|Y\|^2 = d(1+\frac{1}{|Z_0|^2})(1+o(1))$, and $(1+\frac{1}{|Z_0|^2}) = \frac{\|Y\|^2}{d}(1+o(1))$. Substituting in $p(X|Y, Z_0)$, we obtain

$$p(X|Y) \overset{L}{\longrightarrow} N\left(Y\frac{d}{\|Y\|^2}, (1 - \frac{d}{\|Y\|^2})I_d\right). \qquad \square$$

*Heuristic proof for multivariate Cauchy case*

Assume $\frac{\log n}{d^{1/3}} \to 0$ and let $\|Y - X_{[1]}\| \leq \cdots \leq \|Y - X_{[n]}\|$ be the order statistics. We have $w_i \propto (1 + \|Y - X_i\|^2)^{-\frac{d+1}{2}}$, and can write, as usual, $w_{(n)} = \frac{1}{1+T_{n,d}}$, where

$$(17) \quad T_{n,d} = \sum_{j=2}^{n} exp\Big(-\frac{d+1}{2}[\log(1 + \|Y - X_{[j]}\|^2) - \log(1 + \|Y - X_{[1]}\|^2)]\Big).$$

The noise vector follows the multivariate Cauchy distribution, i.e. $\varepsilon = \frac{[Z_1, \ldots, Z_d]^T}{|Z_0|}$, where $Z_0, Z_1, \ldots, Z_d$ are *iid* $N(0, 1)$. Note that

$$1 + \|Y - X_i\|^2 = (1 + \|Y\|^2 + d)\Big(1 + \frac{(\|X_i\|^2 - d) - 2\langle Y, X_i\rangle}{1 + \|Y\|^2 + d}\Big).$$

Now, similarly to the developments in Proposition 3.4, given $Z_0$,

$$1 + \|Y\|^2 + d = d(2 + \frac{1}{Z_0^2})(1 + O_p(d^{-1/2})),$$

by the central limit theorem. Moreover,

$$(18) \quad \log(1 + \|Y - X_i\|^2) - \log(1 + \|Y - X_{[1]}\|^2) = \log(1 + S_i) - \log(1 + S_{(1)}),$$

where

$$\begin{aligned}
S_i &= \frac{\sum_{j=1}^{d}[(X_{ij}^2 - 1) - 2Y_j X_{ij}]}{1 + \|Y\|^2 + d} \\
(19) \quad &= \frac{\sum_{j=1}^{d}[(X_{ij}^2 - 1) - 2Y_j X_{ij}]}{d(2 + \frac{1}{Z_0^2})}(1 + O_p(d^{-1/2})).
\end{aligned}$$

Some calculations show that, conditioning on $Z_0$, the asymptotic variance of $S_i$ is given by $\frac{\sigma^2(Z_0)}{d}$, where $\sigma^2(Z_0) = \frac{6+\frac{4}{Z_0^2}}{(2+\frac{1}{Z_0^2})^2}$. Now, up to a scale factor, the $S_i$'s are the same (conditionally independent) summands as in the Gaussian-Gaussian case. So, we can hope for a uniform Gaussian approximation.

We proceed as if it were exactly Gaussian with mean 0 and variance given by the asymptotic value $\frac{\sigma^2(Z_0)}{d}$. We now expand the logarithm in (19) using $S_i = O_p(d^{-1/2})$. With $W_1, \ldots, W_n$ *iid* $N(0,1)$ and with $W_{(1)}$ the minimum, we obtain

$$exp\big(-\frac{d+1}{2}[\log(1+S_i) - \log(1+S_{(1)})]\big)$$
(20)
$$= exp\{\frac{d^{1/2}}{2}\big([\sigma(Z_0)W_i - \frac{\sigma^2(Z_0)}{2d^{1/2}}W_i^2] - [\sigma(Z_0)W_{(1)} - \frac{\sigma^2(Z_0)}{2d^{1/2}}W_{(1)}^2]\}$$
$$\times (1+o_p(1)).$$

Since $|W_{(1)}| = O_p(\sqrt{\log n})$, we heuristically neglect the cubic term in the expression which is $O_p(|W_{(1)}^3|/d^{1/2}) = o_p(W_{(1)}^2)$.

We now proceed as in the Gaussian-Gaussian case neglecting the $o_p(1)$ term. We have,

$$(21) \quad E(T_{n,d}|W_{(1)}, Z_0) \approx (n-1)q_d^{-1}(W_{(1)})\frac{\int_{W_{(1)}}^\infty q_d(w)\phi(w)\,dw}{\bar\Phi(W_{(1)})},$$

where $q_d(w) = exp\big(-\frac{d^{1/2}}{2}\sigma(Z_0)w + \frac{\sigma^2(Z_0)}{4}w^2\big)$. Letting $A = \frac{\sigma^2(Z_0)}{2}, B = \frac{\sigma(Z_0)}{2}$, using the approximation $\bar\Phi(W_{(1)}) \approx 1$ for the denominator along with Mill's ratio, we compute the integral in (21) to get

$$(22)\quad E(T_{n,d}|W_{(1)}, Z_0) \approx \frac{n-1}{q_d(W_{(1)})(1-A)^{1/2}}\bar\Phi\big[(1-A)^{1/2}W_{(1)}$$
$$+ d^{1/2}B(1-A)^{-1/2}\big]exp\big(\frac{dB^2}{2(1-A)}\big)$$

$$(23)\quad \approx \frac{n-1}{\sqrt{2\pi}(1-A)^{1/2}}exp\big(-\frac{W_{(1)}^2}{2}\big)$$
$$\times \big[(1-A)^{1/2}W_{(1)} + d^{1/2}B(1-A)^{-1/2}\big]^{-1}$$

$$\approx \frac{\sqrt{2\log n}}{\sqrt{2\pi}[(1-A)\sqrt{2\log n} + \sqrt{d}B]}$$

$$(24)\quad = \frac{r_n}{\sqrt{2\pi}[(1-A)r_n + B]}$$

where $r_n = \sqrt{\frac{2\log n}{d}} \to 0$.

We see that collapse occurs (given our approximations) for fixed $Z_0$. We can further calculate the average rate of collapse by taking expectation with respect to $Z_0$. It can be observed that the average rate is dominated by values of $Z_0$ which are near 0. In the case of small $Z_0$, $\sigma(Z_0) = \frac{\sqrt{6+\frac{4}{Z_0^2}}}{2+\frac{1}{Z_0^2}} \approx 2Z_0$ and the normal density is close to $\frac{1}{\sqrt{2\pi}}$. Further $A \approx 1 - 2Z_0^2$ and $B \approx Z_0$. It then follows from (22) that the order of the average rate needed for collapse is $\int_0^\epsilon \frac{r_n}{r_n+z}dz \sim r_n|\log r_n|$, where $\epsilon$ is a small positive number. This is distinctly slower than the Gaussian-Gaussian $r_n$ rate.

These heuristics may be made rigorous if $\log n / d^{1/3} \to 0$. But since the integration by parts remainder terms again dominate, we cannot trace the effect of $Z_0$ precisely.

*Proof of Proposition 4.1.* With $p_Y(y) = \frac{4^{-d/2}}{(2\pi)^{d/2}} exp(-\|y\|^2/4)$, the marginal density of $Y$, write

$$\sum_{i=1}^{n} w_i h(X_i) = \frac{n^{-1}\sum_{i=1}^{n} h(X_i) \frac{exp\left(-\frac{1}{2}\|Y_0 - X_i\|_E^2\right)}{(2\pi)^{d/2} p_Y(Y_0)}}{n^{-1}\sum_{k=1}^{n} \frac{exp\left(-\frac{1}{2}\|Y_0 - X_k\|_E^2\right)}{(2\pi)^{d/2} p_Y(Y_0)}}. \quad (25)$$

Then, the expectation of the numerator under $p(X) \sim N(0, I_d)$,

$$E_0\left[h(X)\frac{exp\left(-\frac{1}{2}\|X - Y_0\|^2\right)}{(2\pi)^{d/2} p_Y(Y_0)}\right] = \int_X h(X) \frac{exp\left(-\frac{1}{2}\|X - Y_0\|^2 - \frac{1}{2}\|X\|^2\right)}{(2\pi)^d p_Y(Y_0)} dX$$

$$= \int_X h(X) \frac{exp\left(-\|X - Y_0/2\|^2\right)}{\pi^{d/2}} dX$$

$$= E_1[h(X|Y_0)].$$

Specializing to $h \equiv 1$, we obtain that the expectation of the denominator in (25) is 1. Now consider the variance of the denominator under $p(X)$,

$$Var_0\left(\frac{1}{n}\sum_{i=1}^{n} \frac{exp\left(-\frac{1}{2}\|X_i - Y_0\|^2\right)}{(2\pi)^{d/2} p_Y(Y_0)}\right) \leq \frac{1}{n}\int_X \frac{exp\left(-\|X - Y_0\|^2\right) exp(-\frac{1}{2}\|X\|^2)}{(2\pi)^{d/2} 4^{-d} exp(-\frac{1}{2}\|Y_0\|^2)} dX$$

$$\leq \frac{1}{n} M^2 (4\sqrt{3})^d.$$

Thus, if $\log n / d \to \infty$, we have $Var_0\left(\frac{1}{n}\sum_{i=1}^{n} \frac{exp(-\frac{1}{2}\|X_i - Y_0\|^2)}{(2\pi)^{d/2} p_Y(Y_0)}\right) \xrightarrow{P} 0$. □

*Proof of Lemma A.2.* By assumption it is sufficient to prove the result for $\lambda_j \equiv 1$, $j = 1, \ldots, d$. We get,

$$E\left|(Z_j + \mu_j)^2 - (1 + \mu_j^2)\right|^k = E\left|(Z_j^2 - 1) + 2\mu_j Z_j\right|^k$$

$$\leq 2^{k-1}\left(E|Z_j^2 - 1|^k + 2^k |\mu_j|^k E|Z_j|^k\right). \quad (26)$$

Let $\lceil x \rceil$ denote the smallest integer greater than or equal to $x$. It is well known that

$$E|Z_j|^k \leq 2^{k/2} \lceil \frac{k}{2} \rceil! \quad (27)$$

Also, with $Z_j'$ an *iid* copy of $Z_j$, then, by Jensen's inequality and (27),

$$E|Z_j^2 - 1|^k \leq E|Z_j^2 - Z_j'^2|^k = 2^k E\left|\frac{Z_j + Z_j'}{\sqrt{2}}\right|^k E\left|\frac{Z_j - Z_j'}{\sqrt{2}}\right|^k = 2^k (E|Z_j|^k)^2$$

$$\leq 4^k (\lceil \frac{k}{2} \rceil!)^2. \quad (28)$$

Now, applying Jensen's inequality again, and noting assumption (8), we have

$$(29) \qquad \sum_{j=1}^{d} |\mu_j|^k \le \Big(\sum_{j=1}^{d} \mu_j^2\Big)^{\frac{k}{2}} d^{-\frac{k-2}{2}} \le \Big(\frac{M}{4}\Big)^{\frac{k}{2}} d.$$

The lemma follows by combining (26) through (29). $\square$

**Acknowledgments.** We are grateful to Chris Snyder, Mesoscale and Microscale Meterology, National Center for Atmospheric Research, for many helpful insights and suggestions.

## References

[1] ANDERSON, J. AND ANDERSON, S. (1999). A monte-carlo implementation of the nonlinear filtering problem to produce ensemble assimilations and forecasts. *Monthly Weather Review* **127** 2741–2758.
[2] ANDERSON, T. (2003). *An Introduction to Multivariate Statistical Analysis*, 3rd ed. Wiley, New York. MR1990662
[3] BENGTSSON, T., SNYDER, C. AND NYCHKA, D. (2003). Toward a nonlinear ensemble filter for high-dimensional systems. *J. Geophysical Research-Atmospheres* **108** 8775.
[4] BERLINER, L. M. (2001). Monte Carlo based ensemble forecasting. *Stat. Comput.* **11** 269–275. MR1842976
[5] BERLINER, M. AND WIKLE, C. (2006). Approximate importance sampling Monte Carlo for data assimilation. In review.
[6] BURGERS, G., P. J., VAN LEEUWEN, P. AND EVENSEN, G. (1998). Analysis scheme in the ensemble Kalman filter. *Monthly Weather Review* **126** 1719–1724.
[7] CRAMÉR, H. (1946). *Mathematical Methods of Statistics*. Princeton Univ. Press.
[8] DONOHO, D. (2000). High-dimensional data analysis: The curses and blessings of dimensionality. *Aide-Memoire of a Lecture at AMS conference on Math Challenges of 21st Centuary.* Available at http://www-stat.stanford.edu/~donoho/Lectures/AMS2000/AMS2000.html.
[9] DOUCET, A., FREITAS, N. AND GORDON, N., EDS. (2001). *Sequential Monte Carlo Methods in Practice*. Springer, New York.
[10] EVENSEN, G. (1994). Sequential data assimilation with a nonlinear quasi-geostrophic model using Monte Carlo methods to forecast error statistics. *J. Geophysical Research* **99** 143–162.
[11] EVENSEN, G. AND VAN LEEUWEN, P. J. (1996). Assimilation of geostat altimeter data for the Agulhas Current using the ensemble Kalman filter. *Monthly Weather Review* **124** 85–96.
[12] FURRER, R. AND BENGTSSON, T. (2007). Estimation of high-dimensional prior and posterior covariance matrices in Kalman filter variants. *J. Multivariate Analysis-Revised* **98** 227–255. MR2301751
[13] GILKS, W. AND BERZUINI, C. (2001). Following a moving target—Monte Carlo inference for dynamic Bayesian models. *J. R. Stat. Soc. Ser. B Stat. Methodol.* **63** 127–146. MR1811995
[14] GORDON, N., SALMON, D. AND SMITH, A. (1993). Novel approach to nonlinear/non-Gaussian Bayesian state estimation. *IEE Proceedings-F* **140** 107–113.

[15] Hamill, T. M., J. S. W. and Snyder, C. (2001). Distance-dependent filtering of background error covariance estimates in an ensemble kalman filter. *Monthly Weather Review* **129** 2776–2790.
[16] Hoar, T., Milliff, R., Nychka, D., Wikle, C. and Berliner, L. (2003). Winds from a Bayesian hierarchical model: Computation for atmosphere-ocean research. *J. Comput. Graph. Statist.* **4** 781–807. MR2037932
[17] Houtekamer, P. and Mitchell, H. (2001). A sequential ensemble Kalman filter for atmospheric data assimilation. *Monthly Weather Review* **129** 123–137.
[18] Houtekamer, P. L. and Mitchell, H. L. (1998). Data assimilation using an ensemble Kalman filter technique. *Monthly Weather Review* **126** 796–811.
[19] Liu, J. (2001). *Monte Carlo Strategies in Scientific Computing*. Springer, New York. MR1842342
[20] Liu, J. and Chen, R. (1998). Sequential Monte Carlo methods for dynamic systems. *J. Amer. Statist. Assoc.* **93** 1032–1044. MR1649198
[21] Molteni, F., Buizza, R., Palmer, T. and Petroliagis, T. (1996). The ECMWF ensemble prediction system: Methodology and validation. *Quarterly J. Roy. Meteorological Society* **122** 73–119.
[22] Pitt, M. and Shepard, N. (1999). Filtering via simulation: Auxilliary particle filters. *J. Amer. Statist. Assoc.* **94** 590–599. MR1702328
[23] Saulis, L. and Statulevicius, V. (2000). *Limit Theorems of Probability Theory*. Springer, New York. MR1798811
[24] Storvik, G. (2002). Particle filters for state-space models with the presence of unknown static parameters. *IEEE Transactions on Signal Processing* **50** 281–289.
[25] Tippett, M. K., Anderson, J. L., Bishop, C. H., Hamill, T. M. and Whitaker, J. S. (2003). Ensemble square-root filters. *Monthly Weather Review* **131** 1485–1490.
[26] Toth, Z. and Kalnay, E. (1997). Ensemble forecasting at NCEP and the breeding method. *Monthly Weather Review* **125** 3297–3319.
[27] van Leeuwen, P. (2003). A variance minimizing filter for large-scale applications. *Monthly Weather Review* **131** 2071–2084.
[28] Wikle, C., Milliff, R., Nychka, D. and Berliner, L. (2001). Spatiotemporal hierarchical bayesian modeling: Tropical ocean surface winds. *J. Amer. Statist. Assoc.* **96** 382–397. MR1939342
[29] Wilson, S. and Stefanou, G. (2005). Bayesian approaches to content-based image retrieval. *Proceedings of the International Workshop/Conference on Bayesian Statistics and Its Applications, Varanasi, India (January 2005)*.